\documentclass[11pt]{article}

\usepackage{amsmath,amsfonts,amssymb}
\usepackage{bm} 
\usepackage{graphicx}
\usepackage{url}
\usepackage{enumerate}

\usepackage[text={14cm,21cm},centering]{geometry}
\renewcommand{\baselinestretch}{1.3}

\usepackage{amsthm}

\theoremstyle{plain}
\newtheorem{theorem}{Theorem}
\newtheorem{lemma}{Lemma}
\newtheorem{cor}[lemma]{Corollary}

\newtheorem{prop}[lemma]{Proposition}
\newtheorem{cond}{Condition}

\theoremstyle{definition}

\newcommand{\indicator}[1]{\mbox{{\rm 1}}\{#1\}}
\newcommand{\Prob}{\mathbb{P}}
\newcommand{\Q}{\mathbb{Q}}
\newcommand{\E}{\mathbb{E}}
\newcommand{\Var}{{\rm Var}}

\title{Asymptotic optimality of
the cross-entropy method for Markov chain problems}
\author{%
Ad Ridder\\
{\it\small Vrije University,}\\
{\it\small Amsterdam, Netherlands}\\
{\it\small aridder@feweb.vu.nl}
}
\begin{document}
\maketitle
%

\begin{abstract}
The correspondence between the cross-entropy method
and the zero-variance approximation to simulate
a rare event problem in Markov chains is shown.
This leads to a sufficient condition that the cross-entropy
estimator is asymptotically optimal.
\end{abstract}

%
\section{Model and problem}\label{s:intro}
We deal with a discrete-time Markov chain
$(X(t))_{t=0}^\infty$ on a finite but large state space $\cal X$,
and with a matrix of transition probabilities $P=(p(x,y))_{x,y\in{\cal X}}$,
i.e., $p(x,y)=\Prob(X(t+1)=y|X(t)=x)$.
The state space is partitioned into three sets:
$\cal G$ is a small set of `good' states, $\cal F$ is a small set of
`failed' or `bad' states, and the other states form the large set
of `internal' states ${\cal T} = {\cal X}\setminus({\cal G}\cup{\cal F})$.
For each internal state
$x\in{\cal T}$ let $\gamma(x)$ be the probability that the Markov chain
will hit the failure set before the good set when the chain
starts in state $x$. Or, more formally, define
\[ T=\inf\{t>0 : X(t)\in{\cal G}\cup{\cal F}\},\]
then $\gamma(x)=\Prob(X(T)\in{\cal F}|X(0)=x)$.
For ease of notation we assume that the good set consists of
a single state, denoted by $\bm{0}$, and that the Markov chain
jumps immediately out of it (but not immediately to a bad state),
i.e., $p(\bm{0},\bm{0})=0$ and $p(\bm{0},{\cal F})=0$.
We are interested in the probability that
when the chain starts in the good state $\bm{0}$,
it will hit the failure set $\mathcal{F}$ before returning to $\bm{0}$.
The associated event is denoted by $A=\indicator{X(T)\in{\cal F}}$.
\par
In this paper we consider the problem of estimating the hitting
probability $\Prob(A)$ by simulation.
The typical applications that we have in mind, are
system failures in models of highly reliable Markovian system,
and excessive backlogs in product-form Jackson queueing networks.
The failure set is a rare event that should occur only with very small probability,
which makes the need for reducing the simulation variance,
for instance by importance sampling.
Reliability and queueing systems have been studied largely in relation to importance
sampling estimation of the performance measures, we refer to overviews
in \cite{heidelberger95} and in \cite{junejashahabuddin06}.
\par
The rare event probability $\Prob(A)$ will be estimated
by an importance sampling simulation method that implements
a change of measure $\Prob^*$, i.e., the estimator $Y$
is the average of i.i.d.\ replications of
\[\indicator{A}\frac{d\,\Prob}{d\,\Prob^*}, \]
where $d\,\Prob/d\,\Prob^*$ is the likelihood ratio,
and where we assume that $\indicator{A}d\,\Prob$ is absolute continuous
w.r.t.\ $d\,\Prob^*$.
The issue is to find a good change of measure in terms
of performance of the estimator $Y$.
The optimal change of measure would be
\[ \Prob^{\rm opt}(\cdot)=\Prob(\cdot|A),\]
for which $Y$ would have zero variance \cite{junejashahabuddin06}.
It is shown in \cite{ecuyertuffin09,glynniglehart89}
that the associated Markov chain
has transition probabilities
\begin{equation}\label{e:opttransprob}
p^{\rm opt}(x,y) = p(x,y)\frac{\gamma(y)}{\gamma(x)},\quad
x\in{\cal T}, y\in{\cal X},
\end{equation}
and $p^{\rm opt}(x,y)=p(x,y)$ for good state $x=\bm{0}$, and bad states
$x\in{\cal F}$. The optimal transition probabilities cannot be
used directly for simulation since they require knowledge of the
unknown hitting probabilities, however they suggest
to construct an importance sampling algorithm
by approximating or estimating the
hitting probabilities $\gamma(x)$.
For instance, let $\Pi(x)$ be the set of all sample paths of the Markov chain
of the form $\pi=(x=x_0\to x_1\to \cdot\to x_k)$ with $x_0,\ldots,x_{k-1}\in{\cal T}$,
$x_k\in{\cal F}$ and $p(x_j,x_{j+1})>0$ ($j=0,\ldots,k-1$), and where $k=1,2,\ldots$.
The probability of path $\pi$ to occur is $p(\pi)=\sum_{j=0}^{k-1}p(x_j,x_{j+1})$,
and, clearly, the hitting probability becomes $\gamma(x) = \sum_{\pi\in\Pi(x)}p(\pi)$.
In \cite{ecuyertuffin09} it is studied how the approximation
\[ \gamma^{\rm app}(x) = \max_{\pi\in\Pi(x)}\,p(\pi) \]
performs in reliability problems. In a slightly different context,
\cite[Section 4]{ecuyeretal10} considers the rare event probability that a random walk
reaches high levels. It is shown that it fits in the Markov chain framework
and an approximation of $\gamma(x)$ is proposed based on an
asymptotic approximation of these probabilities.
\par
Another line of research has been developed in \cite{dupuisetal07,dupuiswang07}
for rare event problems in which the probability of interest
can be approximated via large deviations.
In this framework, the decay rate is given in terms of a variational problem,
or an optimal control problem. These problems are related to a family of nonlinear partial
differential equations known as Hamilton-Jacobi-Bellmann (HJB) equations.
It is shown how subsolutions of the HJB equations associated with rare event problems
could be used to construct efficient importance sampling schemes.
This method has been successfully applied to several queueing systems,
see  \cite{dupuisetal07,dupuiswang07}.
\par
The cross-entropy method for rare event simulation \cite{cebook04}
considers to choose a change of measure $\Prob^{\rm ce}$ from a specified family of
changes of measures that minimizes the Kullback-Leibler distance from
the optimal $\Prob^{\rm opt}$. It has been shown by \cite{deboernicola02}
that the associated transition probabilities $p^{\rm ce}(x,y)$
are of the form
\begin{equation}\label{e:cesol}
p^{\rm ce}(x,y) = \frac{\E[\indicator{A}N(x,y)|X(0)=\bm{0}]}%
{\E\left[\indicator{A}\sum_{z\in{\cal X}}N(x,z)|X(0)=\bm{0}\right]},
\end{equation}
where $N(x,y)$ is the number of times that transition $(x,y)$
occurs in a random sample path of the Markov chain until absorption
in one of the good or bad states.
Again, these transition probabilities cannot be
used directly for simulation since they contain the unknown
variables $N(x,y)$, however, in this case they suggest
to estimate the expectations in expression \eqref{e:cesol},
for instance by simulation. This approach
has been applied to queueing systems in \cite{deboernicola02,deboeretal00},
and to reliability systems in \cite{ridder05}.
\par
In the following sections we shall give more background on the cross-entropy method
and how it results in the expression \eqref{e:cesol}.
More importantly, we shall show that in fact the cross-entropy solution is the
zero-variance distribution, i.e., the matrices of transition
probabilities satisfy $P^{\rm ce}=P^{\rm opt}$.
This identity will be the basis to formulate in Section 3 a sufficient
condition for which the estimated cross-entropy solution
is asymptotically optimal.
\section{Correspondence cross-entropy and zero-variance}\label{s:cezv}
Let $(\Omega,{\cal A},\Prob)$ be the probability space of
the sample paths of the Markov chain $X(0),X(1),\ldots$,
and denote by $\bm{X}$ a random sample path.
Recall that our objective is to execute simulations
of the Markov chain for estimating the rare event probability
$\Prob(A)=\Prob(X(T)\in{\cal F})$, that we execute these simulations
under a change of measure $\Prob^*$,
and that the optimal change of measure is
\begin{equation}\label{e:dpopt}
d\,\Prob^{\rm opt}=\frac{\indicator{A}}{\Prob(A)}d\,\Prob.
\end{equation}
Suppose that we consider only changes of measures $\Prob^*$ under which
the Markov property is retained, say with matrix of transition probabilities
$P^*$, and suppose that we minimize
the Kullback-Leibler distance between these changes of measures and
the optimal one, i.e.,
\[
\inf_{P^* \in {\cal P}}\,{\cal D}(\Prob^{\rm opt},\Prob^*),
\]
where the cross-entropy is defined by
\begin{equation}\label{e:ce}
{\cal D}(\Prob^{\rm opt},\Prob^*)=
\E^{\rm opt}\left[\log\left(\frac{d\Prob^{\rm opt}}{d\Prob^*}(\bm{X})\right)\right]
=\E\left[\frac{d\Prob^{\rm opt}}{d\Prob}(\bm{X})\;
\log\left(\frac{d\Prob^{\rm opt}}{d\Prob^*}(\bm{X})\right)\right].
\end{equation}
Substituting \eqref{e:dpopt}, minimizing the cross-entropy, and deleting constant terms
yields
\begin{equation}\label{e:cesupprogram}
\sup_{P^*\in {\cal P}} \, \E[\indicator{A}\log d\,\Prob^*(\bm{X}) ].
\end{equation}
Since the sample path probability $d\,\Prob^*(\bm{X})$ is a product
of individual transition probabilities, we get
\begin{equation}\label{e:prodform}
d\,\Prob^*(\bm{X}) = \prod_{t=1}^T p^*(X(t-1),X(t))
=\prod_{(x,y)\in{\cal X}\times{\cal X}} p^*(x,y)^{N(x,y)},
\end{equation}
where $N(x,y)$ is the number of times transition $(x,y)$
occurring in the random sample path $\bm{X}$.
Substituting the expression \eqref{e:prodform}
into the cross-entropy optimization program
\eqref{e:cesupprogram}, and applying the first order
condition using a Lagrange multiplier, gives
the solution \eqref{e:cesol} for the individual transition probabilities.
\par
Notice that the optimal transition matrix $P^{\rm opt}$ is a feasible matrix, i.e.,
an element of ${\cal P}$, and thus it must hold that it is the cross-entropy
solution.
We shall give a direct proof of the matrix identity $P^{\rm ce} = P^{\rm opt}$,
based on the expressions of the transition probabilities.
In fact we shall prove a relation between
the expected number of transitions from $x$ to $y$ and
the absorption probability $\gamma(y)$.
Denote by $v(x)$ the expected number of visits to state $x$
starting at $x$ before absorption:
\[
v(x) = \E\left[\sum_{t=0}^\infty\indicator{X(t)=x}\Big|X(0)=x\right].
\]
\begin{prop}\label{p:ngammarelation}
For all $x,y\in{\cal X}$:
\[
\E[\indicator{A}N(x,y)|X(0)=x] = v(x)p(x,y)\gamma(y).
\]
\end{prop}
\begin{proof}
For ease of notation we assume that we have the
equivalent modelling in which all the good and bad states
are absorbing.
Introduce probabilities (for any $x,y\in {\cal X}$)
\begin{align*}
f(x,y)&=\Prob((X(t)) \mbox{ reaches state } y |X(0)=x)\\
g(y)&=\Prob((X(t)) \mbox{ reaches bad set $\cal F$ without a
transition } x\to y|X(0)=y).
\end{align*}
Notice that we allow $x$ and $y$ to be a good or bad state
for which these probabilities are obviously either zero or one.
Consider the event
\[ \{N(x,y)=n\}\,\cap\,\{\mbox{reach bad set from $x$}\},
\]
for $n\geq 1$.
This event can only occur if (A) $n-1$ times
(i) a number of times [transition $x\to y'\neq y$ followed by a return to $x$],
followed by (ii) [transition $x\to y$ followed by a return to $x$];
then (A) is followed by (B) which is
(iii) a number of times [transition $x\to y'\neq y$ followed by a return to $x$],
followed by (iv) [transition $x\to y$ followed by reaching the bad set
without the transition $x\to y$]. That is,
\begin{align*}
\E&[\indicator{A}N(x,y)|X(0)=x]\\
&=
\sum_{n=1}^\infty n
\underbrace{\Bigg(\Big[
\underbrace{\sum_{k=0}^\infty\Big(\sum_{y'\neq y}p(x,y')f(y',x)\Big)^k}_{\rm (i)}
\Big]\underbrace{p(x,y)f(y,x)}_{\rm (ii)}\Bigg)^{n-1}}_{\rm (A)}\\
& \times \;\underbrace{\Big[
\underbrace{\sum_{k=0}^\infty\Big(\sum_{y'\neq y}p(x,y')f(y',x)\Big)^k}_{\rm (iii)}
\Big]\underbrace{p(x,y)g(y)}_{\rm (iv)}}_{\rm (B)}.
\end{align*}
Let us work out the summations using geometric series
and denoting
\[ \alpha = \sum_{y'\neq y}p(x,y')f(y',x);\quad
\beta=\sum_{k=0}^\infty\Big(\sum_{y'\neq y}p(x,y')f(y',x)\Big)^k.
\]
Hence,
\begin{equation}\label{e:nxy}
\E[\indicator{A}N(x,y)|X(0)=x]=\frac{1}{(1-\beta)^2}\,\frac{1}{1-\alpha}
\,p(x,y)g(y).
\end{equation}
In the same manner we determine the absorption probability
\[ \gamma(y)=\Prob((X(t))\mbox{ reaches the bad set }|X(0)=y).\]
Partition this event with respect to the number of transitions
$x\to y$:
\begin{align*}
\gamma(y)&=\sum_{n=0}^\infty \Prob(A; N(x,y)=n|X(0)=y)\\
&=g(y)+\sum_{n=1}^\infty f(y,x)\Prob(A; N(x,y)=n|X(0)=x)\\
&=g(y)+
\sum_{n=1}^\infty
\Bigg(f(y,x)\Big[
\sum_{k=0}^\infty\Big(\sum_{y'\neq y}p(x,y')f(y',x)\Big)^k
\Big]p(x,y)f(y,x)\Bigg)^{n-1}\\
&\times\;\Big[
\sum_{k=0}^\infty\Big(\sum_{y'\neq y}p(x,y')f(y',x)\Big)^k
\Big]p(x,y)g(y)\\
&=g(y)+
\sum_{n=1}^\infty
\Bigg(\Big[
\sum_{k=0}^\infty\Big(\sum_{y'\neq y}p(x,y')f(y',x)\Big)^k
\Big]p(x,y)f(y,x)\Bigg)^{n}
g(y).
\end{align*}
When we include the first term $g(y)$ as the zero-th term of
the summation, we obtain
\begin{equation}\label{e:gammay}
\gamma(y) = \frac{1}{1-\beta}g(y).
\end{equation}
From the expressions \eqref{e:nxy} and \eqref{e:gammay}
we see that
\[
\E[\indicator{A}N(x,y)|X(0)=x] =
\left(\frac{1}{1-\beta}\,\frac{1}{1-\alpha}\right)
p(x,y)\gamma(y).
\]
To conclude, we calculate the proportionality factor:
\begin{align*}
\frac{1}{1-\beta}&\,\frac{1}{1-\alpha} =
\frac{1}{1-\frac{p(x,y)f(y,x)}{1-\alpha}}\,\frac{1}{1-\alpha}\\
&=\frac{1}{1-\alpha-p(x,y)f(y,x)}
=\frac{1}{1-\sum_{y'\neq y}p(x,y')f(y',x)-p(x,y)f(y,x)}\\
&=\frac{1}{1-\sum_{y}p(x,y)f(y,x)}=\frac{1}{1-f(x,x)}=v(x).
\end{align*}
The last equality is a well-known relation for Markov chains,
e.g., see \cite{feller50}.
\end{proof}
\begin{cor}\label{p:ceisopt}
For all $x,y\in{\cal X}$:
\[ p^{\rm ce}(x,y) = p^{\rm opt}(x,y).
\]
\end{cor}
\begin{proof}
The identity follows easily by noting that
\[
\E[\indicator{A}N(x,y)|X(0)=\bm{0}]=f(\bm{0},x)
\E[\indicator{A}N(x,y)|X(0)=x].\]
\end{proof}

\section{Asymptotic optimality}\label{s:ao}
In this section we assume that there is a family of rare events
$\{A_n\}$ parameterized by $n=1,2,\ldots$ such that
each $A_n$ satisfies the model assumptions of the previous section,
and such that $\lim_{n\to\infty}\Prob(A_n)=0$.
Suppose that $Y_n^*$ is an unbiased estimator of
$\Prob(A_n)$ obtained by a change of measure $\Prob^*$.
Then this estimator is asymptotically optimal
if
\[ \lim_{n\to\infty}\,\frac{\log {\E}^*[(Y_n^*)^2]}{\log \Prob(A_n)}
=2,
\]
see for instance \cite{heidelberger95}.
Now, recall the cross-entropy representation $p^{\rm ce}(x,y)$ in
\eqref{e:cesol} of the zero-variance transition probabilities,
and suppose that these are estimated by $\hat{p}^{\rm ce}(x,y)$.
A common approach is to apply an iterative scheme
to the optimization program \eqref{e:cesupprogram}.
Since the program involves the rare event, we first apply a change of measure:
\begin{equation}\label{e:ceprogramcom}
\E[\indicator{X(T)\in{\cal F}}\log d\Prob^*(\bm{X})] =
\E^{(0)}\left[\frac{d\Prob}{d\Prob^{(0)}}
\indicator{X(T)\in{\cal F}}\log d\Prob^*(\bm{X}) \right].
\end{equation}
This is done for a probability measure $\Prob^{(0)}$ such that
(i) the Markov property is retained; (ii) the associated matrix
of transition probabilities is feasible $P^{(0)}\in{\cal P}$;
(iii) the set ${\cal F}$ is `not so' rare under $\Prob^{(0)}$.
The program \eqref{e:ceprogramcom} is solved iteratively by estimation:
let $\bm{X}^{(1)}, \ldots, \bm{X}^{(k)}$ be i.i.d.\ sample paths
of the Markov chain generated by simulating the states according
to a matrix of transition probabilities $P^{(j)}$ until absorption
in the good or bad states, then we calculate for $j=0,1,\ldots$
\begin{equation}\label{e:ceupdate}
P^{(j+1)} = \arg\max_{P^*\in {\cal P}} \;\frac{1}{k}\sum_{i=1}^{k}
\frac{d\Prob}{d\Prob^{(j)}}(\bm{X}^{(i)})
\indicator{X^{(i)}(T)\in{\cal F}}\log d\Prob^*(\bm{X}^{(i)}).
\end{equation}
We repeat this `updating' of the change of measure a few times
until the difference $P^{(j+1)}-P^{(j)}$ is small
enough (in some matrix norm).
For details we refer to \cite{cebook04}.
\par
In this way we obtain an implementable
change of measure $\hat{\Prob}^{\rm ce}$, and its associated importance
sampling estimator is denoted by
\begin{equation}\label{e:ceestimator}
\hat{Y}^{\rm ce}_n = \frac{d\Prob}{d\hat{\Prob}^{\rm ce}}(\bm{X})\,
\indicator{A_n}.
\end{equation}
Clearly, this estimator is unbiased, i.e.,
$\hat{\E}^{\rm ce}[\hat{Y}^{\rm ce}_n]=\Prob(A_n)$\footnote{We denote
the expectation w.r.t.\ measure $\Prob$ by $\E$,
w.r.t.\ measure $\hat{\Prob}^{\rm ce}$ by
$\hat{\E}^{\rm ce}$, w.r.t.\ measure $\Prob^{\rm opt}$ by $\E^{\rm opt}$, etc.}.
We claim that it is asymptotically optimal if the following condition holds.
\begin{cond}\label{c:aocond}
There are finite posivite constants $K_1,K_2$ such that for all $n$
\[
K_1\leq {\cal D}(\Prob^{\rm opt},\hat{\Prob}^{\rm ce})\leq K_2.
\]
\end{cond}
\noindent
This is a condition on the approximation of the zero-variance
measure by the implementation of cross-entropy method.
Actually, a weaker condition for the
upper bound suffices:
${\cal D}(\Prob^{\rm opt},\hat{\Prob}^{\rm ce})=o(\log\Prob(A_n))$
for $n\to\infty$.
\begin{theorem}
Assume Condition \ref{c:aocond}. Then the cross-entropy importance sampling estimator
\eqref{e:ceestimator} is asymptotically optimal.
\end{theorem}
\begin{proof}
Notice that ${\cal D}(\Prob^{\rm opt},\hat{\Prob}^{\rm ce})
=\E^{\rm opt}[\log d\Prob^{\rm opt}/d\hat{\Prob}^{\rm ce}(\bm{X})]\geq 0$.
Because $\log\Prob(A_n)\to -\infty$, the upper bound in
Condition \ref{c:aocond} ensures that
\[
\lim_{n\to\infty}\frac{\E^{\rm opt} \left[ \log
\frac{d\Prob^{\rm opt}}{d\hat{\Prob}^{\rm ce}}(\bm{X})\right]}%
{\log \Prob(A_n)}=0.
\]
Furthermore, the lower bound gives
\[
\limsup_{n\to\infty}
\frac{ \log \E^{\rm opt} \left[
\frac{d\Prob^{\rm opt}}{d\hat{\Prob}^{\rm ce}}(\bm{X})\right]}%
{\E^{\rm opt} \left[ \log
\frac{d\Prob^{\rm opt}}{d\hat{\Prob}^{\rm ce}}(\bm{X})\right]}
<\infty.
\]
Now consider, (using \eqref{e:dpopt} in the third equality in the following lines),
\begin{align*}
\hat{\E}^{\rm ce}\left[(\hat{Y}^{\rm ce}_n)^2\right]&=
\hat{\E}^{\rm ce}\left[\left(\frac{d\Prob}{d\hat{\Prob}^{\rm ce}}(\bm{X})\,
\indicator{A_n}\right)^2\right]\\
&=
\hat{\E}^{\rm ce}\left[\left(\frac{d\Prob}{d\Prob^{\rm opt}}(\bm{X})\,
\indicator{A_n}\right)^2\,
\left(\frac{d\Prob^{\rm opt}}{d\hat{\Prob}^{\rm ce}}(\bm{X})\right)^2\right]\\
&=
\Prob(A_n)^2\, \hat{\E}^{\rm ce}\left[
\left(\frac{d\Prob^{\rm opt}}{d\hat{\Prob}^{\rm ce}}(\bm{X})\right)^2\right]
=
\Prob(A_n)^2\, \E^{\rm opt}\left[
\frac{d\Prob^{\rm opt}}{d\hat{\Prob}^{\rm ce}}(\bm{X})\right].
\end{align*}
So, we can conclude
\begin{align*}
& \frac{\log \hat{\E}^{\rm ce}[(\hat{Y}^{\rm ce}_n)^2]}{\log \Prob(A_n)}
=\frac{\log (\Prob(A_n))^2 + \log \E^{\rm opt}\left[
\frac{d\Prob^{\rm opt}}{d\hat{\Prob}^{\rm ce}}(\bm{X})\right]}%
{\log \Prob(A_n)}\\
&= 2 + \frac{\log \E^{\rm opt}\left[
\frac{d\Prob^{\rm opt}}{d\hat{\Prob}^{\rm ce}}(\bm{X})\right]}%
{\log \Prob(A_n)},
\end{align*}
with
\[
\lim_{n\to\infty}
\frac{\log \E^{\rm opt} \left[
\frac{d\Prob^{\rm opt}}{d\hat{\Prob}^{\rm ce}}(\bm{X})\right]}%
{\log \Prob(A_n)}=
\lim_{n\to\infty}
\frac{\log \E^{\rm opt} \left[
\frac{d\Prob^{\rm opt}}{d\hat{\Prob}^{\rm ce}}(\bm{X})\right]}%
{\E^{\rm opt} \left[ \log
\frac{d\Prob^{\rm opt}}{d\hat{\Prob}^{\rm ce}}(\bm{X})\right]}%
\;
\frac{\E^{\rm opt} \left[ \log
\frac{d\Prob^{\rm opt}}{d\hat{\Prob}^{\rm ce}}(\bm{X})\right]}%
{\log \Prob(A_n)} = 0.
\]
\end{proof}
\section{A numerical example}\label{s:ex}
We illustrate the theorem by the simple example of simulating the
$M/M/1$ queue (Poisson-$\lambda$ arrivals, exponential-$\mu$ services)
where $\lambda<\mu$.
We consider its associated discrete-time Markov chain $(X(t))_{t=0}^\infty$
by embedding the continuous-time queueing process at the jump times.
The transition probabilities are
\begin{align*}
p &= p(x,x+1)=\frac{\lambda}{\lambda+\mu}\quad (x=0,1,\ldots)\\
q &= p(x,x-1)=\frac{\mu}{\lambda+\mu}\quad (x=1,2,\ldots).
\end{align*}
The rare event is hitting state $n$ before returning to the zero state.
For this model the optimal (zero-variance) transition
probabilities follow easily from calculating the hitting probabilities
\[ \gamma(x) = \Prob((X(t)) \mbox{ reaches $n$ before $0$}| X(0)=x),\]
for $x=1,\ldots,n-1$, by solving the equations
\[ \gamma(x) = p(x,x-1)\gamma(x-1) + p(x,x+1)\gamma(x+1),\]
with boundary conditions $\gamma(0)=0$ and $\gamma(n)=1$.
Let $\sigma=\mu/\lambda$. Then we get
\[
\gamma(x)=\frac{1-\sigma^x}{1-\sigma^n};\quad
p^{\rm opt}(x,x+1)=p\frac{1-\sigma^{x+1}}{1-\sigma^x};
\quad
p^{\rm opt}(x,x-1)=q\frac{1-\sigma^{x-1}}{1-\sigma^{x}}.
\]
Notice that $p^{\rm opt}(1,2)=1$.
\par
The cross-entropy  between the optimal probability measure
$\Prob^{\rm opt}$ and any other probability measure $\Q$
which is associated with transition probabilities
$q(x,x+1)$ and $q(x,x-1)$, can be determined as follows
(where we apply the product form \eqref{e:prodform}):
\begin{align*}
{\cal D}& (\Prob^{\rm opt},\Q)
=\E^{\rm opt}\left[\log \frac{d\Prob^{\rm opt}}{d\Q}(\bm{X})\right]\\
&=\E^{\rm opt}\left[\prod_{(x,y)\in{\cal X}\times{\cal X}} \log
\left(\frac{p^{\rm opt}(x,y)}{q(x,y)}\right)^{N(x,y)} \Big| X(0)=0\right]\\
&=\sum_{(x,y)\in{\cal X}\times{\cal X}} %
\log \frac{p^{\rm opt}(x,y)}{q(x,y)} \,
\E^{\rm opt}[N(x,y) | X(0)=0 ].
\end{align*}
We may follow the same reasoning as in the proofs of Propositions
\ref{p:ngammarelation} and \ref{p:ceisopt} for calculating
$\E^{\rm opt}[N(x,y) | X(0)=0 ]$. Notice that under $\Prob^{\rm opt}$
$\indicator{A}=1$, $f(0,x)=1$, and $\gamma(y)=1$, thus
\[
\E^{\rm opt} [N(x,y) | X(0)=0 ]
=f^{\rm opt}(0,x)v^{\rm opt}(x)p^{\rm opt}(x,y)\gamma^{\rm opt}(y)
=v^{\rm opt}(x)p^{\rm opt}(x,y).
\]
Finally, the visiting numbers $v^{\rm opt}(x)=1/(1-f^{\rm opt}(x,x))$
can be calculated numerically via recursion (using $f^{\rm opt}(x,x+1)=1$):
\begin{align*}
f^{\rm opt}(x,x)&=p^{\rm opt}(x,x-1) + p^{\rm opt}(x,x+1)f^{\rm opt}(x+1,x)\\
f^{\rm opt}(x+1,x)&=
\frac{p^{\rm opt}(x+1,x)}{1- p^{\rm opt}(x+1,x+2)f^{\rm opt}(x+2,x+1)}.
\end{align*}
\par
We have implemented the cross-entropy method of Section \ref{s:ao} for queueing
parameters $\lambda=0.8$ and $\mu=1$. The rare-event state $n$ was increased from
$n=10$ until $n=250$. The cross-entropy updating rule \eqref{e:ceupdate}
was iterated ten times, starting from the uniform transition probabilities
(the sample sizes $k$ were increased proportionally to $n$).
After the ten updating iterations we estimated the rare event probability
$\Prob(A_n)$ with sample size 1000 and collected the estimated
relative errors (RE)
$\sqrt{\widehat{\Var}^{\rm ce}[\hat{Y}^{\rm ce}_n]}/\hat{\E}^{\rm ce}[\hat{Y}^{\rm ce}_n]$
of the associated estimators, and the estimated ratios (RAT)
$\log \hat{\E}^{\rm ce}[(\hat{Y}^{\rm ce}_n)^2]/\log \hat{\E}^{\rm ce}[\hat{Y}^{\rm ce}_n]$.
The figures below illustrate that
${\cal D}(\Prob^{\rm opt},\hat{\Prob}^{\rm ce})/\Prob(A_n)\to 0$ as $n\to\infty$,
and that RAT is close to two, meaning that the estimator is
asymptotically optimal. The last figure with the relative errors RE
indicates even strong efficiency (bounded RE).
\par\bigskip\noindent
\begin{center}
\includegraphics[width=8cm, height=5cm]{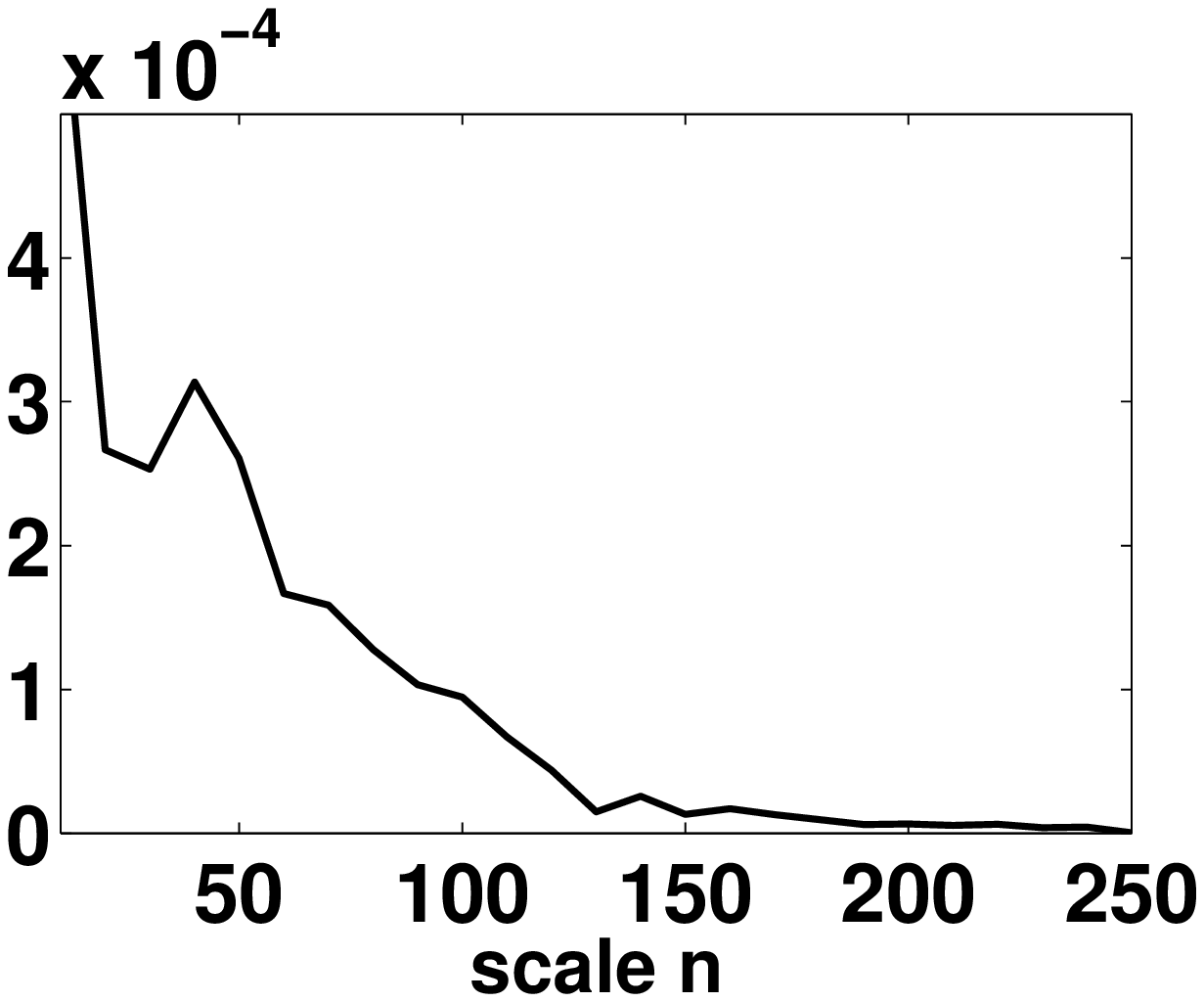}
\par\smallskip
{\it Figure 1. ${\cal D}(\Prob^{\rm opt},\hat{\Prob}^{\rm ce})/\Prob(A_n)$.}

\bigskip

\includegraphics[width=8cm, height=5cm]{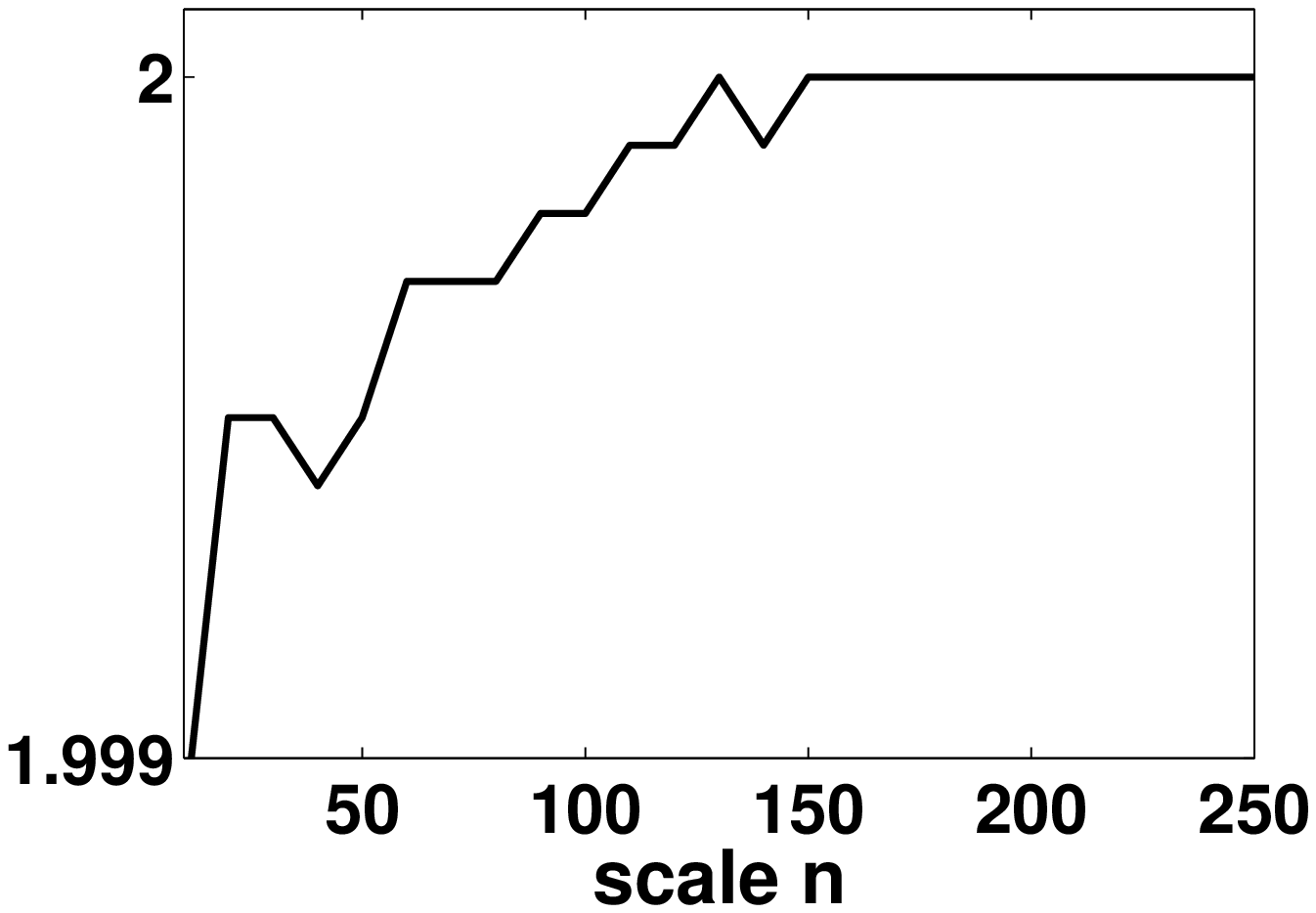}
\par\smallskip
{\it Figure 2. Estimated ratio RAT.}

\bigskip

\includegraphics[width=8cm, height=5cm]{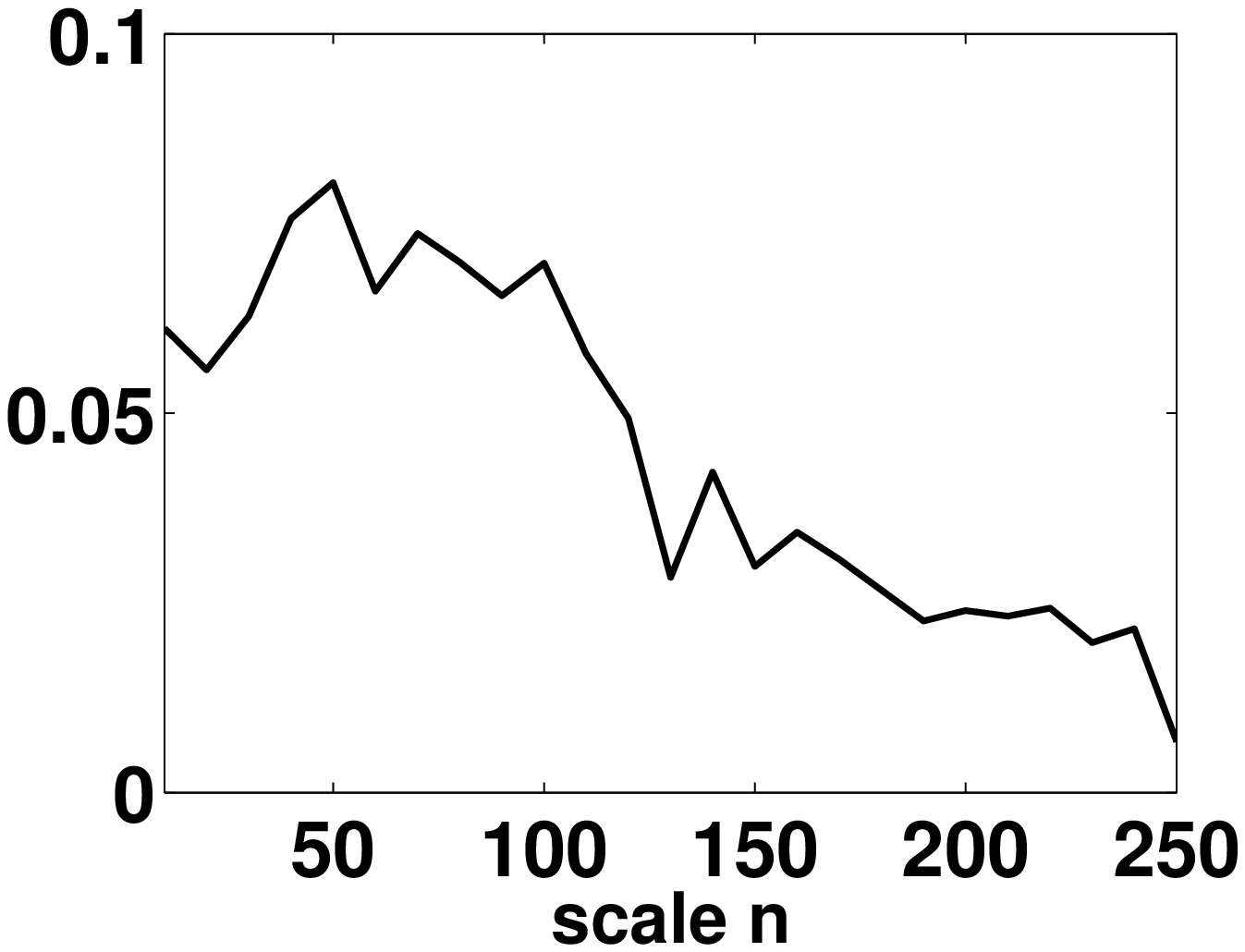}
\par\smallskip
{\it Figure 3. Estimated relative error RE.}
\end{center}
\section{Conclusion}\label{s:conclusion}
After having constructed an importance sampling algorithm
for rare-event simulation, a key issue is to assess the
statistical performance of the associated importance sampling
estimator. In this paper we have considered rare-event problems
in Markov chains, for which we have used the cross-entropy method
as the engine of finding a change of measure for executing the
importance sampling simulations.
We have shown that for these problems the cross-entropy method
coincides with the zero-variance approach. This non-implementable optimal change
of measure is estimated by an implementable change of measure
that is returned by the cross-entropy method.
Our main result is that we give a sufficient condition for the
associated importance sampling estimator to be logarithmically efficient.
Further investigations are undertaken to obtain conditions for strong
efficiency.
\section*{Acknowledgements}
The author would like to thank Bruno Tuffin for his helpful discussions,
and his hospitality during a visit to INRIA Rennes Bretagne Atlantique.

\large
\renewcommand{\baselinestretch}{1.0}

\end{document}